\pdfoutput=1
\documentclass{article}
\usepackage[margin=25mm]{geometry}
\usepackage{amsmath}
\usepackage{amsfonts}
\usepackage{amssymb}
\usepackage{graphicx}
\usepackage[utf8]{inputenc}
\usepackage[T2A]{fontenc}
\usepackage[ukrainian,russian,english]{babel}
\usepackage[unicode, pdftex]{hyperref}
\usepackage{xcolor}
\usepackage[all]{xy}

\newtheorem{theorem}{Theorem}

\title{Structure of Morse flows with at most six singular points on the torus with a hole}

\author{Maria Loseva, Alexandr Prishlyak, Kateryna Semenovych, Dariia Synieok}

\begin{document}
\newcounter{contnumfig}
\setcounter{contnumfig}{0}

\maketitle
\begin{abstract}
   We describe all possible topological structures of Morse flows and typical gradient saddle-nod bifurcation of flows  on the 2-dimensional torus with a hole in the case that the number of singular point of flows is at most six. To describe structures, we use separatrix diagrams of flows. The saddle-node bifurcation is specified by selecting a separatrix in the separatrix diagram of the flow befor the bifurcation.
\end{abstract}
\section*{Introduction}

In this paper, we consider Morse flows on the torus with a hole. It is a gradient flow of Morse function. Since the function increases along each trajectory, the flow has no cycles and polycycles. In general position, a typical gradient flow is a Morse flow (Morse-Smale flow without closed trajectories). In typical one-parameter families of gradient flows, two types of bifurcations are possible: saddle-node and saddle connection. The corresponding vector fields at the time of the bifurcation are fields of codimension 1. In our case, they completely determine the topological type of the bifurcation.
To classify Morse flows, a separatrix diagram is often used, in which separetrices are trajectories that belonging to one-dimensiona stable or unstable manifolds.  We apply this approach to the classification of typical saddle-nod bifurcation.

Without loss of generality, we assume that under bifurcation (as the parameter increases), the number of singular points does not increase. The saddle-node bifurcation is defined by a separatrix, which is contructed to a point. We mark this separatrix on the diagram. 

We colar  stable separatrices in red  and unstable separatrices in green.

 Kronrod \cite{Kronrod1950} and Reeb \cite{Reeb1946} construct topological invariants of functions oriented 2-maniofolds. It was generlized in \cite{lychak2009morse} for  non-orientable two-dimensional manifolds and in   \cite{Bolsinov2004, hladysh2017topology, hladysh2019simple} for manifolds with boundary, in \cite{prishlyak2002morse} for non-compact manifolds. 

If we consider Morse flows as gradient flows of Morse functions and fix the value of functions in singular points, the flow determinate the topological structure of the function \cite{lychak2009morse, Smale1961}. Therefore, Morse--Smale flows structure is closely related to the structure of the functions.

Possible structures of smooth function on closed 2-manifolds was considered in \cite{bilun2023morseRP2, hladysh2019simple, hladysh2017topology,  prishlyak2002morse, prishlyak2000conjugacy,  prishlyak2007classification, lychak2009morse, prishlyak2002ms, prish2015top, prish1998sopr,  bilun2002closed,  Sharko1993, prish2002Met}, on 2-manifolds with the boundary in \cite{hladysh2016functions, hladysh2019simple, hladysh2020deformations} and on closed 3- and 4-manifolds in  \cite{prishlyak1999equivalence, prishlyak2001conjugacy, hlp2023}.

In \cite{bilun2023gradient, Kybalko2018, Oshemkov1998, Peixoto1973, prishlyak1997graphs, prishlyak2020three, akchurin2022three, prishlyak2022topological, prishlyak2017morse,  kkp2013,  prish2002vek,  prishlyak2021flows,  prishlyak2020topology,   prishlyak2019optimal, prishlyak2022Boy}, 
the structures  of flows on closed 2- manifolds and 
\cite{bilun2023discrete, loseva2016topology, prishlyak2017morse, prishlyak2022topological, prishlyak2003sum, prishlyak2003topological, prishlyak1997graphs, prishlyak2019optimal} on manifolds with the boundary were described.
Topological properties of Morse-Smale flows on 3-manifolds was investigated in \cite{pbp2023, prish1998vek,  prish2001top, Prishlyak2002beh2, prishlyak2002ms,   prish2002vek, prishlyak2005complete, prishlyak2007complete, hatamian2020heegaard, bilun2022morse, bilun2022visualization}.


The purpose of this paper is to describe all possible topological structures of the Morse flows   with no more than six singular points on the torus with a hole.



 
\section{Morse-Smale flow and separatrix diagram}

Typical vector fields on compact 2-manifolds are Morse-Smale fields. Morse fields or Morse-Smale gradient-like fields that do not contain closed trajectories are tipical among the gradient fields. They satisfy three properties:

1) singular points are nondegenerate;

2) there are no separatric connections in the interior of the manifold;

3) $\alpha$-limiting ($\omega$-limiting) set of each trajectory is a singular point.

\begin{figure}[ht!]
\center{\includegraphics[width=0.95\linewidth]{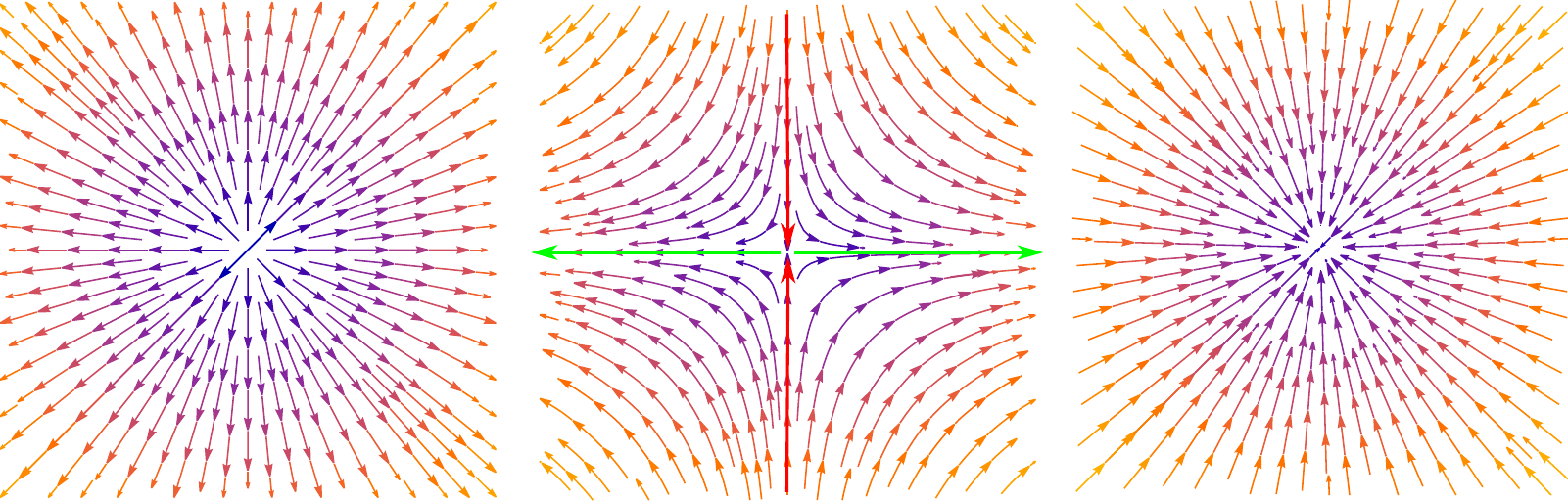}
}
\put (-380,-12) {1}
\put (-230,-12) {2}
\put (-70,-12) {3}
\caption{Interior singular points: 1) a source, 2) a saddle,  3) a sink.}
\label{nsn}
\end{figure}

There are tree type of singular points for Morse flow on 2-manifolds: 1) a source, 2) a saddle,  3) a sink. In the neighbourhood of a source vector field is topologically equivalent to the field $X=\{x,y\}$, in saddle to $X=\{x, -y\}$  and in sink to $X=\{-x, -y\}$. \textit{Separatrix} is a trajectories, that go in or go out from the saddle. \textit{Separatrix diagram } of Morse flow  is the 2-manifold with following elements in them: 1) singular points, 2) separatrices, 3) boundary trajectories.  

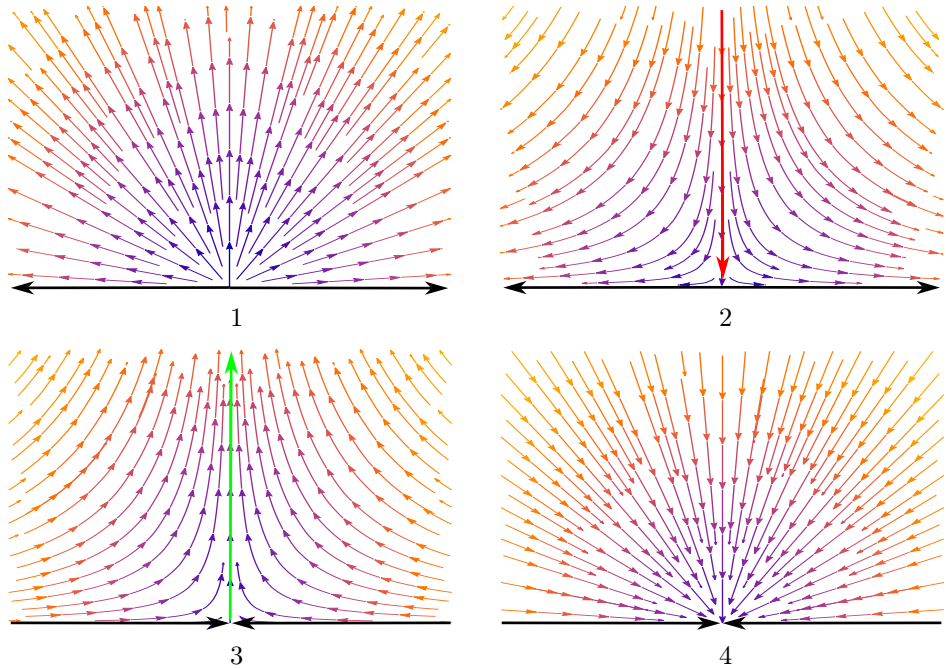
\begin{figure}[ht!]
\center{\includegraphics[width=0.75\linewidth]{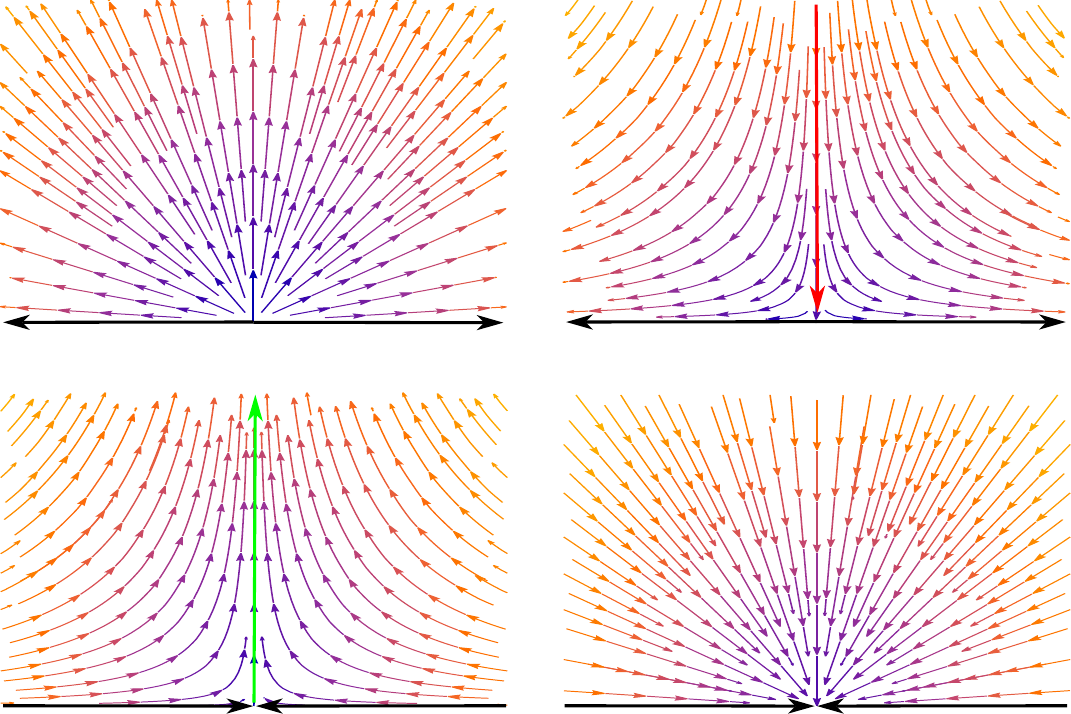}
}
\put (-270,115) {1}
\put (-85,115) {2}
\put (-270,-13) {3}
\put (-85,-13) {4}
\caption{Boundary singular points: 1) a source, 2) a saddle, 3) a saddle, 4) a sink.}
\label{bnssn}
\end{figure}

We colour separatrices in the interior of 2-manifold in two colour: red if $\omega$-limit set of the separatrix is a saddle and in green if $\alpha$-limit set of the separatrix is a saddle.

\section{The structure of Morse flows on the torus with the hole}

In this section, we consider the structure of Morse flows with at most six singular points on a torus with a hole. When restricting the Morse flow to the boundary, we obtain a Morse flow on the boundary that is homeomorphic to the circle. Two types of points are possible for Morse flow on a circle:  sources and sinks. Since the number of sources is equal to the number of sinks, the total number of points on the boundary is even. So it is equal to 2, 4 or 6.

When doubling the manifold, the interior points are doubled, but the border points are not. For example, if there were three internal singular points, and two on the boundary, then there will be $3\times 2+ 2=8$ singular points on the doubling. At the same time, all these points are non-degenerate. If $N$ is the number of sources and sinks, and $S$ is the number of saddles on the doubling, which is a surface of genus 2, then according to the theorem
Poincaré-Hopf $S-N=2$. 

Since each flow has at least one source and at least one sink, then $N \geq 2$ and $S+N \geq 6$.
\subsection{Flows with 4 singular points}
 Consider flows with 4 singular points on a torus with a hole. If they all lie on the boundary, then the redouble has the same number of singular points, which contradicts the fact that the doubling has at least 6 special points.

So a flow with 4 singular points on a torus with a hole has two points on the boundary and two points inside. The source and sink lie on the boundary, otherwise the condition $S-N=2$ for doubling is violated.

\begin{figure}[ht!]
\center{\includegraphics[width=0.95\linewidth]{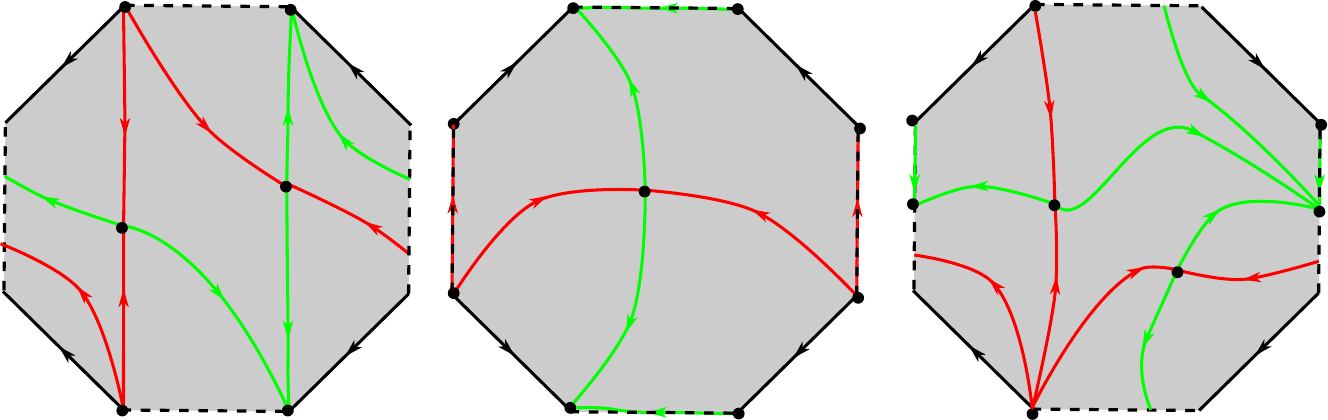}
}
\put (-380,-10) {1}
\put (-230,-10) {2}
\put (-70,-10) {3}
\caption{Morse flows  with 4 and 5 singular points}
\label{t-45}
\end{figure}

We will show that there is a single flow whose separatrix diagram is shown in Fig. 1.1 (here, octagons have the upper sides glued to the lower ones, and the left sides to the right ones). To do this, we will cover the hole with a 2-disc and fill the 2-disc with trajectories starting at the source and ending at the sink. We get a polar Morse flow on the torus. Since all such flows are topologically equivalent to [10], and all regular trajectories corresponding to the 2-disc are equivalent, we have that all Morse flows on a torus with a hole with 4 singular points have a diagram equivalent to the diagram in Fig. 3.1.

\subsection{Flows with 5 singular points}

According to the reasoning above, two types of flows with 5 singular points are possible: 1) 4 points lie on the boundary and one point is interior, 2) two points lie on the boundary and three are interior.

In the first case, the source and the sink lie on the boundary. We  draw separatrix from the saddles on the boundary. This can be done in a single way with accuracy up to homeomorphism (topological equivalence). Then the inner saddle is inside the octagon, and it is possible to draw separatrices from it to the sinks and sources on the boundary by only one way. Therefore, all flows of the first type are topologically equivalent. Their diagram is shown in fig. 3.2.

In the second case, the flow has one source and one sink, and one of the saddles lies on the boundary. Let, for certainty, the source lies on the boundary, and the sink is inside. Let's draw a green trajectory from the saddle on the boundary to the sink, which can be done uniquely, with accuracy up to homeomorphism. By compressing this trajectory to a point, we get a flow with 4 singular points, which is unique. The inverse procedure, transforming the sink on the boundary to the saddle on the boundary, the sink inside, and the green separatrix, can also be carried out uniquely. Therefore, such a single flow is possible. It is shown in fig. 3.3.

If we change the directions of all trajectories, we will get an inverted flow with a sink on the boundary and an internal source.

Thus, on a torus with a hole there exist three topologically non-equivalent flows with 5 singular points.

\subsection{Flows with 6 singular points}
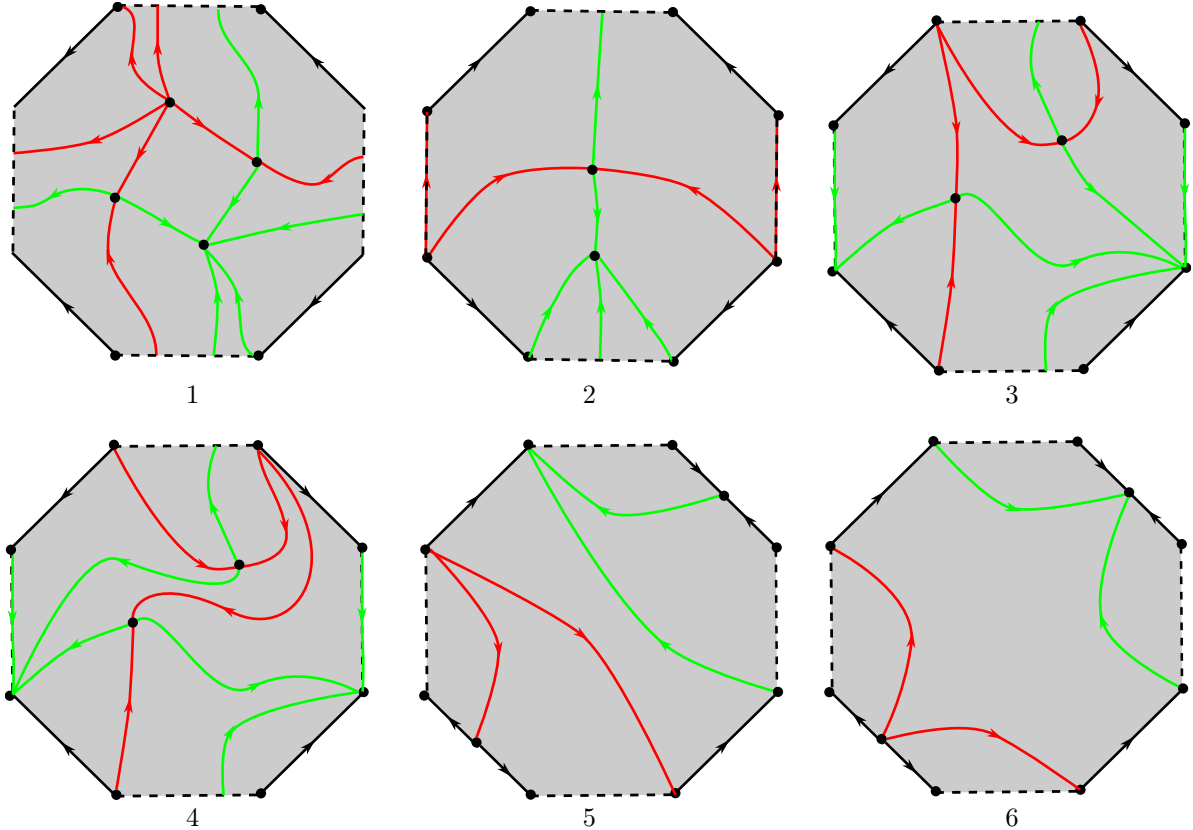
\begin{figure}[ht!] 
\center{\includegraphics[width=0.95\linewidth]{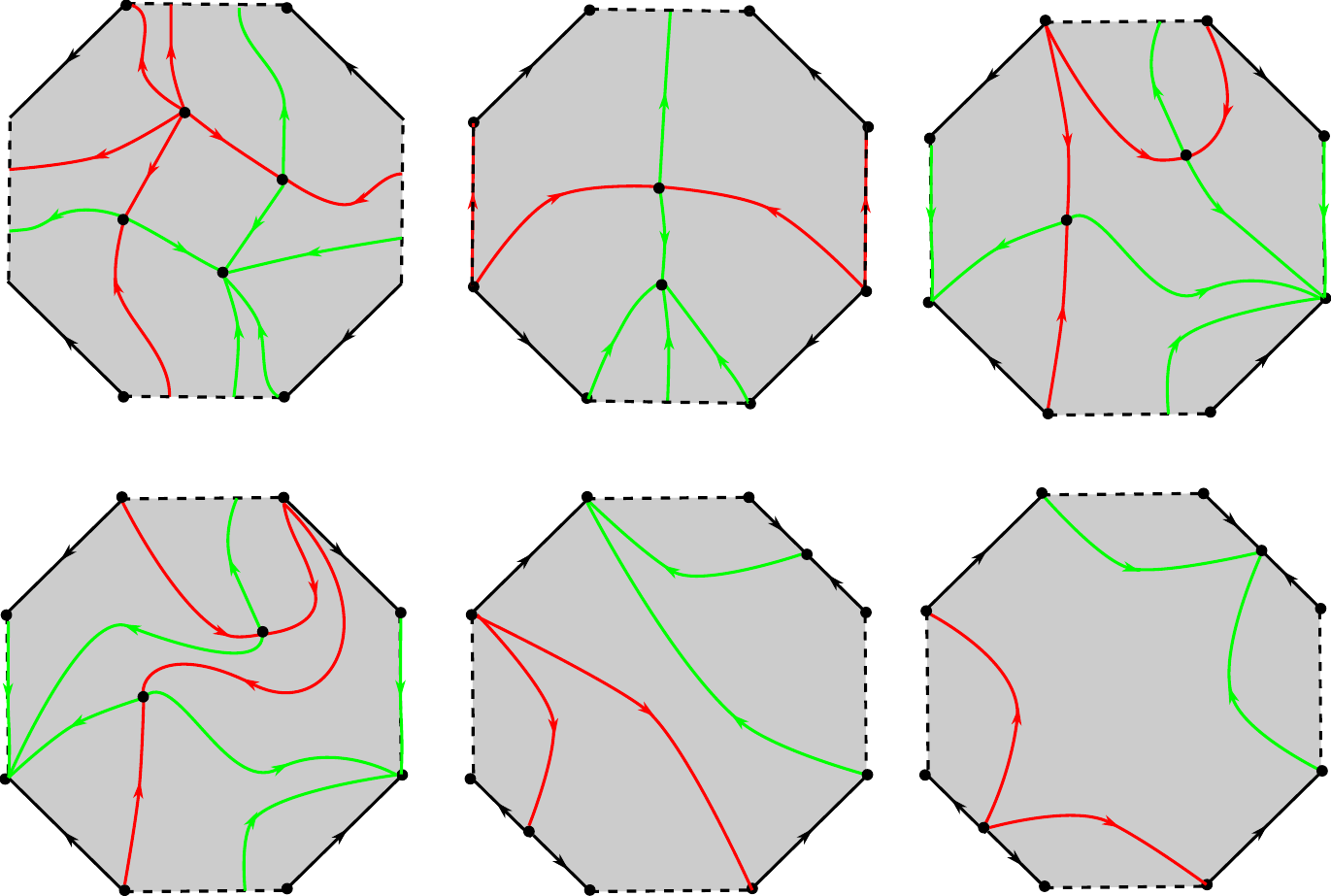}
}
\put (-380,150) {1}
\put (-230,150) {2}
\put (-70,150) {3}
\put (-380,-10) {4}
\put (-230,-10) {5}
\put (-70,-10) {6}
\caption{Morse flows with 6 singular point (part 1)}
\label{T2-ms6}
\end{figure}

If the Morse flow has 6 singular points, then the following types of them are possible:

1) All 6 points belong to the boundary. Then one of them is a sink, one a source, and 4 saddles.

2) Four points belong to the boundary, and two to the interior. Internal points are saddle points.

3) Four points belong to the boundary, and two to the interior. Among the interior points, only one is a saddle.

4) Two points lie on the boundary and both of these points are saddle points.

5) Two points lie on the boundary, one of them is the source, and the other is the sink. Among the inner points, there are three saddles, and one is a source or sink.

Consider the flows of the first type. Two cases are possible for them: a) there is a trajectory at the boundary between the source and the sink (Fig. 4.5), such a trajectory does not exist (Fig. 4.6). In each of these cases, the separatrices are carried out uniquely with accuracy up to the homeomorphism.

For flows of the second type, there is only one saddle point on the boundary. Among the other three points on the boundary, there are two sources and one sink or two sinks and one source. We will consider only the first option (two sources). The second option will be obtained with the first by rotation of the directions on all trajectories. 5 green separatrix come out of the source. A separatrix that enters the saddle at the boundary breaks the rest of the separatrix in the source into two parts. If there are 2 separators in each of them, then the only such flow is shown in fig. 4.4, if there are three separatrixes in one particle, and one in the other, then the separatrix diagram is shown in fig. 4.3. The case when all separatrix belong to one part is impossible , because otherwise stable manifolds of saddle points will intersect, which is impossible.

For flows of the third type, we will assume that a non-saddle internal point is a sink.  Then  one source and three saddles lie on the boundary. All green separatrices are carried out unambiguously, and therefore the red ones are also unambiguous. Therefore, we have one flow with the diagram in Fig. 4.2.

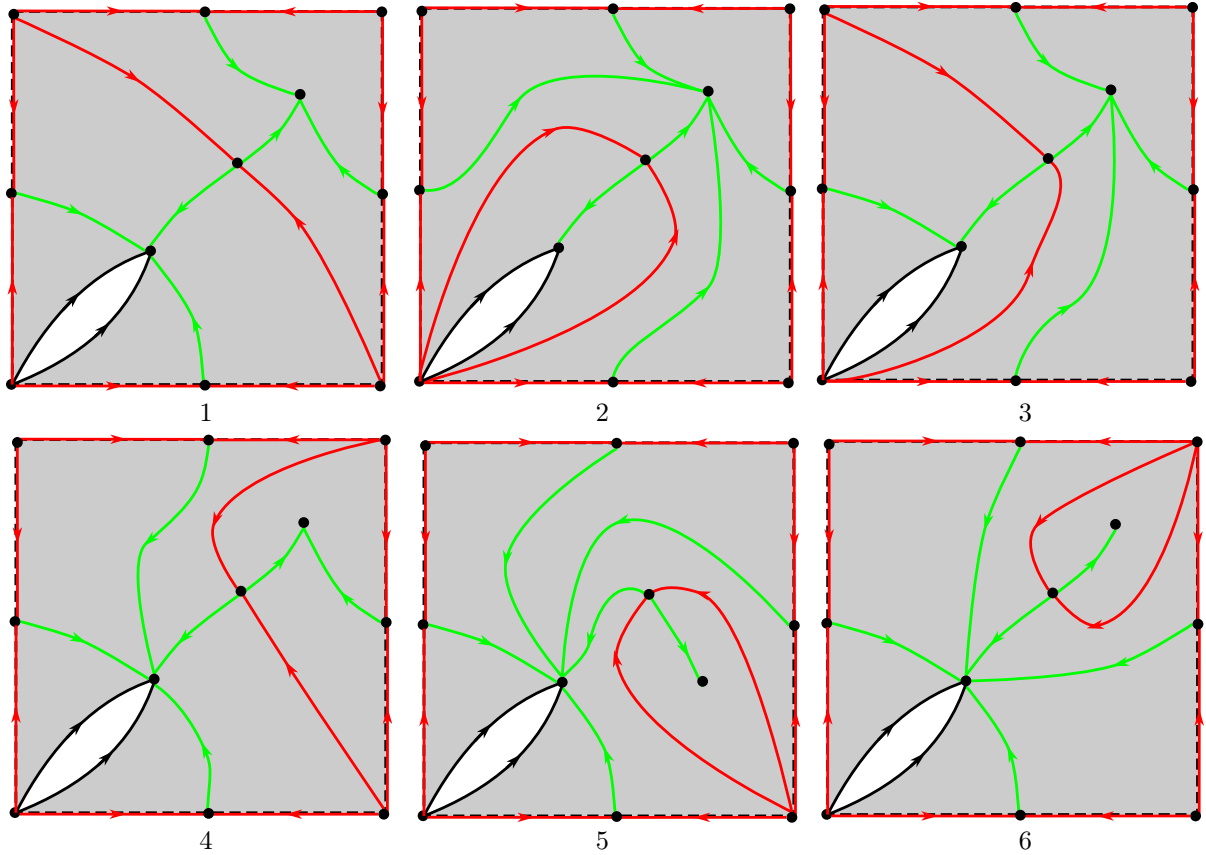
\begin{figure}[ht!]
\center{\includegraphics[width=0.96\linewidth]{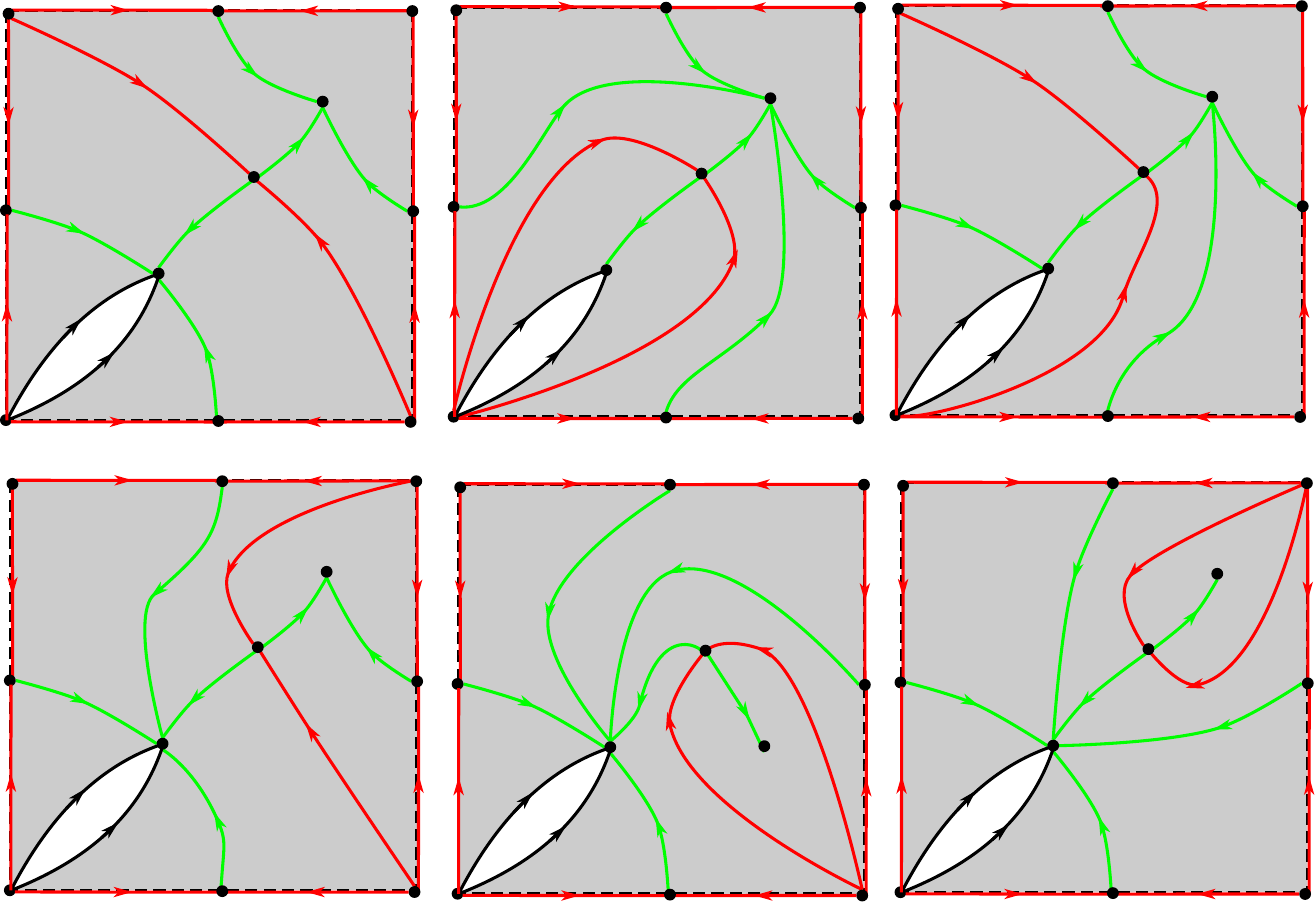}
}
\put (-380,152) {1}
\put (-230,152) {2}
\put (-70,152) {3}
\put (-380,-10) {4}
\put (-230,-10) {5}
\put (-70,-10) {6}
\caption{Morse flows with 6 singular points (part 2)}
\label{T-ms6a}
\end{figure}

For flows of the fourth type, there are two saddles , a source and a sink inside. This flow can be obtained from a polar Morse flow on a torus by cutting an arc of a regular trajectory. The ends of the cut will be two saddle points. Since all such arcs are the same (one can be mapped into another using topological equivalence), we can get only one flow of this type. It is shown in fig. 4.1.

By analogy with the previous one, flows of type 5 will be obtained from Morse flows on a torus with two sources, one drain and three saddles, cutting the torus along one of the regular trajectories between the source and the drain. A total of 6 variants of such cuts are possible. They are shown in fig. 5.

We went through all the possible options, and therefore it is fair

\begin{theorem}
The following possible structures of Morse flow exist on the 2-dimensional torus with a hole:
\begin{itemize}
\item
with four singular points -- one flow (fig. 3.1);
\item
with five singular points -- tree flows (fig.3.2, fig. 3.3 and reverse fig. 3.3);
\item
with six singular points  -- 21 flows (fig. 4 and fig 5).
\end{itemize}

y
\end{theorem}

\section{Saddle-nod bifurcations and boundary saddle bifurcations}

In this section, we consider all possible gradient bifurcations in which the separatrix or trajectory on the boundary contracts to a point. In order for the gradient field to be formed, other separatrices or boundary trajectories must not form loops. This means that only those trajectories that are not multiple edges on the separatrix diagram can participate in the bifurcation. In particular, if the boundary contains only two singular points, then the trajectories of the boundary do not specify bifurcations. The following types of bifurcations are possible:

SN -- separatrix, between internal points, one of which is a saddle,

BSN -- trajectory on the boundary, between points, one of which is a saddle,

HN -- separatrix between the saddle on the boundary and the internal point,

HS -- separatrix between the inner saddle and a point on the boundary,

BDS -- a trajectory on the border, between two saddles.

Next, we provide a list of bifurcations for each diagram.

Fig. 3.1: No possible bifurcation.

Fig. 3.2: 2 BSN, BDS.

Fig. 3.3: HN. (x2)

Fig. 4.1: 2HN.

Fig. 4.2: BSN, BDS, HN. (x2)

Fig. 4.3: 2HS, 2BSN, BDS. (x2)

Fig. 4.4: 4HS, 2BSN, BDS. (x2)

Fig. 4.5:  2BSN, 3 BDS. 

Fig. 4.6:  2 BSN,  BDS. 

Fig. 5.1: 2 SN, HS.(x2)

Fig. 5.2:  SN, HS.(x2)

Fig. 5.3: 3 SN,  HS.(x2)

Fig. 5.4: 2 SN, 2 HS.(x2)

Fig. 5.5:  SN,  HS.(x2)

Fig. 5.6:  SN, HS.(x2)

Summing up the corresponding numbers, we have the values in the table. 1.

\section*{Conclusion}

All possible structures of Morse flows and typical one-parameter gradient saddle-nod bifurcations on the torus with a hole in which no more than six singular points are found (see Table 1). We hope that the research carried out in this paper can be extended on flows with a larger number of singular points and on other surfaces.

\begin{table}[ht]
	\centering
		\begin{tabular} {|c|c|c|c|c|c|c|c|c|c|}
		\hline
Number of points
& 
Morse & SN &  BSN & BDS & HN & HS  
\\ 

\hline
4& 1 & 0 & 0 & 0 & 0 & 0 
 \\
\hline
5 & 3 & 0 & 2  & 1 & 2 & 0    
 \\
\hline
6& 21 & 20& 14 & 10 & 4 & 26 
 \\

 \hline		

		\end{tabular}
	\caption{Number of Morse flows and bifurcations on the torus with a hole }
	\label{tab:NF}
\end{table}


\textsc{Taras Shevchenko National University of Kyiv}

Maria Loseva \ \ \ \ \ \ \  \ \ \ \ \ \ \ \
\textit{Email:} \text{ mv.loseva@gmail.com} \   \ \ \   \ \ \ 
\textit{ Orcid ID:} \text{0000-0002-2282-206X}

Alexandr Prishlyak \ \   \ \ \  
\textit{Email address:} \text{ prishlyak@knu.ua} \ \ \ \ 
\textit{ Orcid ID:} \text{0000-0002-7164-807X}

Kateryna Semenovych \ \ 
\textit{Email:} \text{ kateryna.semenovych@knu.ua}  \ 
\textit{ Orcid ID:} \text{0000-0001-9717-1524}

Dariia Synieok \ \ \ \ \ \ \ \ \ \ \ 
\textit{Email:} \text{ dasha.sineok@knu.ua}  \   \ \ \  \ \ \ \ \
\textit{ Orcid ID:} \text{0009-0002-4594-370X}
\end{document}